\documentclass{article}
\usepackage{amsmath}
\newcommand{\openbox}{\leavevmode
  \hbox to.77778em{%
  \hfil\vrule
  \vbox to.675em{\hrule width.6em\vfil\hrule}%
  \vrule\hfil}}
\newenvironment{proof}[1][Proof]
{\par\addvspace{6pt}\normalfont \itshape #1\@{.}\hskip\labelsep\ignorespaces\normalfont}
{\hfill\openbox\par\addvspace{6pt}}
\def\b{\backslash}
\def\tr{\textrm}

\def\A{\mathcal{A}}
\def\B{\mathcal{B}}
\def\mc{\mathcal}
\def\sep{[n]^{(r)}_{*}}
\def\sepk{[n]^{(r)}_{k}}

\newtheorem{theorem}{Theorem}

\newtheorem{lemma}{Lemma}                    
\newtheorem{conjecture}{Conjecture}

\begin{document}
\title{Intersecting Families of Separated Sets
}
\author{John Talbot
\\Merton College
\\University of Oxford
\\E-mail: talbot@maths.ox.ac.uk}
\date{March 10th 2002}
\maketitle
\begin{abstract}  
A set $A\subseteq \{1,2,\ldots,n\}$ is said to be \emph{$k$-separated} if, when considered on the circle, any two elements of $A$ are separated by a gap of size at least $k$.

We prove a conjecture due to Holroyd and Johnson \cite{HOL},\cite{HOL2} that an analogue of the Erd\H os-Ko-Rado theorem holds for $k$-separated sets. In particular the result holds for the vertex-critical subgraph of the Kneser graph identified by Schrijver \cite{SCH}, the collection of separated sets. We also give a version of the Erd\H os-Ko-Rado theorem for weighted $k$-separated sets.
\end{abstract}
\section{Introduction}
A family of sets is \emph{intersecting} if any two sets from the
family meet. If $[n]^{(r)}$ is the collection of all $r$-sets from
$[n]=\{1,2,\dots ,n\}$ how large can an intersecting family of sets
$\mc{A}\subseteq [n]^{(r)}$ be? If $n < 2r$ this is easy to answer
since $[n]^{(r)}$ is intersecting. However, for $n \geq 2r$ this
question is more difficult. It was answered by Erd\H os, Ko and Rado
\cite{EKR}.
\begin{theorem}[Erd\H os, Ko and Rado \cite{EKR}]
\label{EKRthm}
Let $n\geq 2r$ and $\mathcal{A} \subset [n]^{(r)}$ be intersecting. Then $|\mathcal{A}| \leq |\mathcal{A}_{1}|$, where $\mathcal{A}_{1}=\{A\in [n]^{(r)}:1\in A\}$.
\end{theorem} 
{\footnotesize Mathematics Subject Classification: 05D05 Extremal Set Theory, 05C65 Hypergraphs.}
\newpage
When $n\geq 2r+1$ the above result can be extended to show that equality holds iff $\A\simeq \A_1$.

Many questions concerning families of sets from $[n]^{(r)}$ can be framed in the language of graphs. Consider the \emph{Kneser graph}, $K_{n,r}$, with vertex set $[n]^{(r)}$ and edges between any two vertices corresponding to disjoint $r$-sets. The Erd\H os-Ko-Rado theorem (Theorem \ref{EKRthm}) can be restated as: the largest independent set of vertices in $K_{n,r}$ has order $\binom{n-1}{r-1}$.

One of the most fundamental properties of a graph is its chromatic number. A
longstanding conjecture due to Kneser \cite{KNE} was that the chromatic number
of $K_{n,r}$ is $n-2r+2$. This was answered in the affirmative by Lov\'asz in
1977 \cite{LOV}. Later Schrijver \cite{SCH} identified a vertex-critical
subgraph of $K_{n,r}$, that is a minimal subgraph of $K_{n,r}$ with chromatic number $n-2r+2$. In order to describe this subgraph we require the following definition. We say that a set $A \in [n]^{(r)}$ is \emph{separated} if, when considered as a subset of $[n]$ arranged around a circle in the usual ordering, $A$ does not contain any two adjacent points. Schrijver's vertex-critical subgraph of the Kneser graph is the subgraph induced by those vertices corresponding to the collection of all separated sets in $[n]^{(r)}$.

Let us denote the collection of all separated sets in $[n]^{(r)}$ by $\sep$. Then the corresponding subgraph of the Kneser graph has (by Schrijver's result) chromatic number $n-2r+2$. However, the independence number of this subgraph was previously not known. This was a rather curious situation since generally determining the independence number of a graph is ``easier'' than determining its chromatic number. 

For the remainder of this paper we will consider a well-known conjecture of Holroyd and Johnson on this problem, namely that an analogue of the Erd\H os-Ko-Rado theorem holds for intersecting families of separated sets. In fact their conjecture is more general. They define for any integer $k\geq 1$ the collection of $k$-\emph{separated} $r$-sets in $[n]^{(r)}$ to be those $r$-sets $A=\{a_{1}, \dots ,a_{r}\}$ with $a_1<a_2 < \cdots <a_r$ satisfying $a_{i+1}-a_{i}>k$ for $i=1, \dots ,r$, where $a_{r+1}=a_{1}+n$. We will denote this family by $[n]^{(r)}_{k}$. Note that a $1$-separated set is simply a separated set.
\begin{conjecture}[Holroyd and Johnson \cite{HOL},\cite{HOL2}]
\label{HOLcon}
\emph{Let $n$, $k$ and $r$ be positive integers satisfying $n\geq (k+1)r$. Suppose $\mathcal{A} \subseteq \sepk $ is intersecting. Then $|\mathcal{A}| \leq |\mathcal{A}_{1}^{\ast}|$, where $\A_{1}^{*}=\{A\in \sepk: 1\in A\}$.}
\end{conjecture}
Our main result, Theorem \ref{kCproof}, is a proof of this conjecture. The key idea used in the proof is a type of ``compression''. Standard proofs using compression (such as the original proof of the Erd\H os-Ko-Rado theorem \cite{EKR}) generally rely on certain properties of sets being preserved under this operation. Our proof is quite different since the property of being $k$-separated is not preserved under our compression operation. 

In the next section of this paper we prove the conjecture for separated sets (the case k=1). This proof is then easily generalised in the subsequent section to yield the full result. We also characterise the extremal families. 

In the final section of this paper we give a new version of the Erd\H os-Ko-Rado theorem for weighted $k$-separated sets.
\section{The Erd\H os-Ko-Rado Theorem for separated sets}
\begin{theorem}
\label{Cproof}
Let $n \geq 2r$ and $\mc{A}\subseteq \sep$ be intersecting then $|\A|\leq |\A_1^*|$ where $\mc{A}_1^*=\{A\in\sep: 1 \in A \}$. Moreover, for $n\neq 2r+2$ the extremal family $\A_1^*$ is unique up to isomorphism. 

If $n=2r+2$ then other extremal families exist. For example if $d=\lfloor r/2\rfloor$ then for each $1\leq i \leq d$ the following family is extremal $$\mc{B}_i=\{A\in \sep:|A\cap\{1,3,\ldots,4i+1\}|\geq  i+1\}.$$
\end{theorem}

Before giving the proof in full we sketch the basic ideas used. 

The proof uses induction on $n$, where $\mc{A}\subseteq\sep$ is our intersecting family. The base case ($n=2r$) is trivial so we may assume that $n\geq 2r+1$. We show that $|\mc{A}|=|\mc{A}_1|+|\mc{A}_2|$, where $\mc{A}_1 \subseteq [n-1]_*^{(r)}$ and $\mc{A}_2\subseteq [n-2]_*^{(r-1)}$ are both intersecting. Since $n\geq 2r+1$ the inductive hypothesis applies in both cases and this yields the desired result.

In order to describe the families $\mc{A}_1$ and $\mc{A}_2$ we introduce a ``compression'' function, $f$, that maps points in $[n]\b\{1\}$ anticlockwise by one and leaves $1$ fixed. Then $\mc{A}_1$ consists of those sets in $f(\mc{A})$ that are still separated (as subsets of $[n-1]$). Clearly $\mc{A}_1$ is an intersecting family in $[n-1]_*^{(r)}$.

If $A\in\sep$ but $f(A)\not\in [n-1]_*^{(r)}$ then $A$ must either contain the pair $\{1,3\}$ or the pair $\{2,n\}$. So let $\mc{D}$ consist of those sets in $\mc{A}$ containing $\{1,3\}$ and let $\mc{E}$ consist of those sets in $\mc{A}$ containing $\{2,n\}$. Applying the compression to both of these families we see that each set in $f(\mc{D})$ contains $\{1,2\}$ and each set in $f(\mc{E})$ contains $\{1,n-1\}$. We now remove the point 1 from each of these compressed sets to give a new family $\mc{G}=(f(\mc{D})-\{1\})\cup(f(\mc{E})-\{1\})$. Clearly this is a disjoint union, the fact that $\mc{G}$ is also intersecting is less obvious.

We now need to consider those cases when two distinct sets $A,B\in\mc{A}$ both compress to the same set (i.e. $f(A)=f(B)$). It is easy to see that this can only happen if $A\Delta B=\{1,2\}$. Hence $|\mc{D}|=|f(\mc{D})|$ and $|\mc{E}|=|f(\mc{E})|$. Also any set in $\mc{A}_1$ has at most two preimages in $\mc{A}$. Let $\mc{F}$ consist of those sets in $\mc{A}_1$ that have two preimages in $\mc{A}$ then each set in $\mc{F}$ contains $1$. Consider the new family  $\mc{H}=(\mc{F}-\{1\})\cup\mc{G}$. Then $\mc{H}$ is the union of disjoint families (to see this just check that if $H\in \mc{H}$ then $H\cap\{2,n-1\}=\emptyset \iff H \in \mc{F}-\{1\}$). Furthermore $\mc{H}$ is also intersecting. Now no set in $\mc{H}$ contains $1$ and so $|\mc{H}|=|f(\mc{H})|$. Applying our compression function to $\mc{H}$ we define $\mc{A}_2=f(\mc{H})$. We then find that $\mc{A}_2$ is in fact an intersecting family in $[n-2]_*^{(r-1)}$. 

Finally $|\mc{A}|=|\mc{A}_1|+|\mc{F}|+|\mc{D}|+|\mc{E}|=|\mc{A}_1|+|\mc{F}|+|\mc{G}|=|\mc{A}_1|+|\mc{H}|=|\mc{A}_1|+|\mc{A}_2|$ as claimed.

In order to show the uniqueness of the extremal family (for $n\neq 2r+2$) we note first that the result holds trivially for $n=2r$ and $n=2r+1$ so we may assume that $n \geq 2r+3$. We then proceed by induction on $r$. Clearly the result holds for $r=1$ so we proceed to the inductive step. Consider the family $\A_2$ defined above. For equality to hold we must have $|\A_2|=\binom{n-r-2}{r-2}$ and since $n-2 \geq 2(r-1)+3$ our inductive hypothesis for $r-1$ implies that there exists $i\in [n]$ such that every set in $\A_2$ contains $i$. Considering the different possibilities for $i$ we find that either $\A\simeq \A_1^*$ or there exists $j \in [n]$ such that $\{A\in\sep:j,j+2\in A \}\subset \mc{A}$. Without loss of generality we may suppose that $j=1$. This easily implies that every set in $\mc{A}$ must contain either $1$ or $3$. Partition $\A$ as $\mc{A}=\B_1\cup\B_3\cup\B_{1,3}$, where $\B_1=\{A\in\A:1\in A,3\not\in A\}$, $\B_3=\{A\in\A:1\not\in A,3\in A\}$ and $\B_{1,3}=\{A\in\A:1,3\in A\}$ . To prove that $\A\simeq \A_1^*$ it is sufficient to show that one of the families $\B_1$ or $\B_3$ must be empty. It can be shown that if $n\geq 2r+3$ and both families are non-empty then $\A$ contains two disjoint sets, a contradiction.
\begin{proof}[Proof of Theorem \ref{Cproof}]
We proceed by induction on $n$. First note that the result is trivial for $n\leq 4$ so we may suppose $n\geq 5$. Also we may suppose that $n \geq 2r+1$ since if $n=2r$ then there are only two separated sets and these are disjoint so the result follows.

Let $\mc{A}\subseteq \sep$ be intersecting. Consider the following partition of $\A$
\[
\mc{A}=\mc{B}\cup\mc{C}\cup\mc{D}\cup\mc{E},
\]
where
\begin{eqnarray*}
\mc{B}& =&\{A\in \A: 1 \not\in A\tr{ and }(2\not\in A\tr{ or }n\not\in A)\},\\
\mc{C}& =&\{A\in \A: 1 \in A\tr{ and }3\not\in A\},\\
\mc{D}& =&\{A\in \A: 1,3\in A\},\\
\mc{E}& =&\{A\in \A: 2,n\in A\}.
\end{eqnarray*}
Define the function $f:[n]\to [n-1]$ by 
\[
f(j)=\left\{\begin{array}{cc} 1,& j=1\\ j-1,& j\geq 2.\end{array}\right.
\]
We will need to use the following results (their proofs are given later).
\begin{lemma}
\label{cval}
If $\mc{A}_1^*=\{A\in\sep:1\in A\}$ then \[|\mc{A}_1^*|=\binom{n-r-1}{r-1 }.\]
\end{lemma}
\begin{lemma}
\label{fequal}
If $A,B\in\sep$, $A\neq B$ and $f$ is as defined above then $$f(A)=f(B) \implies A\Delta B=\{1,2\}.$$
\end{lemma}
\begin{lemma}
\label{fsep}
If $A\in \mc{B}\cup\mc{C}$ and $\mc{B},\mc{C},f$ are as defined above then $f(A)\in [n-1]_*^{(r)}$.
\end{lemma}
\begin{lemma}
\label{sint}
If $\mc{I}$ is an intersecting family of sets then so is $f(\mc{I})$.
\end{lemma}
Since no set in $\mc{B}$ contains $1$ and no set in $\mc{C}$ contains $2$, Lemma \ref{fequal} implies that $|f(\mc{B})|=|\mc{B}|$ and $|f(\mc{C})|=|\mc{C}|$. Hence 
\[
|\mc{B}|+|\mc{C}|=|f(\mc{B})|+|f(\mc{C})|=|f(\mc{B})\cup f(\mc{C})|+|f(\mc{B})\cap f(\mc{C})|.
\]
Then, since $\mc{B}\cup\mc{C}$ is an intersecting family in $\sep$, Lemmas \ref{fsep} and \ref{sint} imply that $f(\mc{B})\cup f(\mc{C})=f(\mc{B}\cup\mc{C})$ is an intersecting family in $[n-1]_*^{(r)}$. Hence by our inductive hypothesis for $n-1 \geq 2r$ and Lemma \ref{cval} we have 
\[
|f(\mc{B})\cup f(\mc{C})|\leq \binom{n-1-r-1}{r-1}=\binom{n-r-2}{r-1}.\]
Let \[
\mc{F}=f(\mc{B})\cap f(\mc{C}).\]
For any family of sets $\mc{G}$ let \[
\mc{G}-\{1\}=\{G\b\{1\}:G\in\mc{G}\}.\]
Define \[
\mc{H}=(f(\mc{D})-\{1\})\cup(f(\mc{E})-\{1\})\cup(\mc{F}-\{1\}).\]
We now require the following result.
\begin{lemma}
\label{slem}
Let $\mc{D},\mc{E},\mc{F},\mc{H}$ and $f$ be as defined above, then 
\begin{itemize}
\item[(a)]$\mc{H}\subseteq [n-1]^{(r-1)}_*$.
\item[(b)]$f(\mc{D})-\{1\}$, $f(\mc{E})-\{1\}$ and $\mc{F}-\{1\}$ are pairwise disjoint families of sets.
\item[(c)] $\mc{H}$ is intersecting. 
\item[(d)] $f(\mc{H}) \subseteq [n-2]_*^{(r-1)}$. 
\end{itemize}
\end{lemma}
Since no set in $\mc{H}$ contains $1$ Lemma \ref{fequal} implies that $|\mc{H}|=|f(\mc{H})|$. Lemma \ref{slem}(c) says that $\mc{H}$ is intersecting and so Lemma \ref{sint} implies that $f(\mc{H})$ is also intersecting. Moreover Lemma \ref{slem}(d) tells us that $f(\mc{H}) \subseteq [n-2]_*^{(r-1)}$. So using Lemma \ref{cval} and the inductive hypothesis for $n-2\geq2(r-1)$ we obtain\[
|\mc{H}|\leq \binom{n-2-(r-1)-1}{(r-1)-1}=\binom{n-r-2}{r-2}.\]
Using Lemma \ref{fequal} again we also have $|\mc{D}|=|f(\mc{D})|$ and $|\mc{E}|=|f(\mc{E})|$. Finally Lemma \ref{slem}(b) tells us that
\begin{eqnarray*}
|\mc{H}|&=&|f(\mc{D})|+|f(\mc{E})|+|\mc{F}|\\
 &=&|\mc{D}|+|\mc{E}| + |\mc{F}|.
\end{eqnarray*}
Hence
\begin{eqnarray}
\nonumber
|\mc{A}|&=&|\mc{B}|+|\mc{C}|+|\mc{D}|+|\mc{E}|\\
\nonumber
& = & |f(\mc{B})\cup f(\mc{C})|+|\mc{F}|+|\mc{D}|+|\mc{E}|\\
\nonumber
& = & |f(\mc{B})\cup f(\mc{C})|+|\mc{H}|\\
\label{bound}
& \leq & \binom{n-r-2}{r-1} + \binom{n-r-2}{r-2}=\binom{n-r-1}{r-1}=|\mc{A}_1^*|.
\end{eqnarray}
This completes the proof of the bound in Theorem \ref{Cproof}. 

In order to prove that the extremal family is unique for $n\neq 2r+2$ we will require the following results (again their proofs follow later).
\begin{lemma}
\label{uni1}
If equality holds in (\ref{bound}) and the family $f(\mc{H})$ defined above is isomorphic to $\B_1^*=\{A\in [n-2]_*^{(r-1)}:1\in A\}$ then either $\A\simeq \A_1^*$ or there exists $j\in[n]$ such that $\A_{j,j+2}^*=\{A\in\sep:j,j+2\in A \}\subset \A$.
\end{lemma}
\begin{lemma}
\label{uni2}
If equality holds in (\ref{bound}) and $A_{1,3}^*=\{A\in \sep:1,3\in A\} \subset \A$ then every set in $\A$ meets $\{1,3\}$. Moreover if $A=\{1,a_2,\ldots,a_r\}\in\sep$ satisfies $A\cap\{1,3\}=\{1\}$ then exactly one of the sets $A$ and $g(A)=\{3,a_2+1,\ldots,a_r+1\}$ belongs to $\A$. 
\end{lemma}
Suppose equality holds in (\ref{bound}). If $n=2r$ or $n=2r+1$ it is easy to deduce that $\A\simeq \A_1^*$ so we may suppose that $n\geq 2r+3$. We will prove that $\A\simeq \A_1^*$ by induction on $r$. Clearly the result holds for $r=1$ or $2$ so suppose $r\geq 3$. Since equality holds in (\ref{bound}) we must have $|f(\mc{H})|=\binom{n-r-2}{r-2}$ and so our inductive hypothesis for $r$ implies that $f(\mc{H})\simeq \B_1^*=\{A\in [n-2]_*^{(r-1)}:1\in A\}$. If $\A\not\simeq \A_1^*$ then Lemma \ref{uni1} implies that there exists $j\in [n]$ such that $\A^*_{j,j+2}=\{A\in\sep:j,j+2\in A\}\subset \A$. Without loss of generality we may suppose that $j=1$.

By the first part of Lemma \ref{uni2} we know that every set in $\A$ meets $\{1,3\}$. Hence we may write $\mc{A}=\B_1\cup\B_3\cup\B_{1,3}$, where $\B_1=\{A\in\A:1\in A, 3 \not\in A\}$, $\B_3=\{A\in\A:3\in A, 1 \not\in A\}$ and $\B_{1,3}=\{A\in\A:1,3\in A\}$. For uniqueness to hold we need to show that one of the families $\B_1$ or $\B_3$ must be empty. We will prove that if both families are non-empty then $\A$ contains two disjoint sets, a contradiction.

Since we are supposing that equality holds in (\ref{bound}), Lemma \ref{uni2} tells us that if $A=\{1,a_2,\ldots,a_r\}\in\sep$ and $A\cap\{1,3\}=\{1\}$ then exactly one of the sets $A$ and $g(A)=\{3,a_2+1,\ldots,a_r+1\}$ belong to $\A$. We now wish to show that if $A\in \B_1$ and $B\in\sep$ is obtained from $A$ by shifting one vertex of $A\b\{1\}$ one place clockwise then $B\in \B_1$. Suppose not, then by Lemma \ref{uni2} we have $g(B)\in\B_3$. But $g(B)\cap A =\emptyset$, contradicting the fact that $\A$ is intersecting. Similarly if $A\in \B_3$ and $B\in\sep$ is obtained from $A$ by shifting one vertex of $A\b \{3\}$ one place anti-clockwise then $B\in \B_3$.

Now suppose, for a contradiction, that both $\B_1$ and $\B_3$ are non-empty then by shifting points one by one we may suppose that $A=\{3,5,\ldots,2r+1\}\in\B_3$ and $B=\{1,n-2r+3,\ldots,n-3,n-1\}\in \B_1$. If $n$ is odd these are disjoint and we have a contradiction. If $n$ is even then either there is a set in $\B_1$ not containing $n-1$ and so $C=\{1,n-2r+2,\ldots,n-4,n-2\}\in\B_1$ and $A\cap C=\emptyset$ or, assuming $2r+1\leq n-2$, $D=\{1,5,7,\ldots,2r+1\}\not\in \A$. If the latter holds then Lemma \ref{uni2} implies that $g(D)=\{3,6,8,\ldots,2r+2\}\in\B_3$. But $B\cap g(D)=\emptyset$. This contradiction completes the proof of Theorem \ref{Cproof}. \end{proof}
We now prove the lemmas. The definitions of $\mc{A},\mc{B},\mc{C},\mc{D},\mc{E},\mc{F},\mc{H}$ and $f$ are as in the proof of Theorem \ref{Cproof}.
\begin{proof}[Proof of Lemma \ref{cval}]
Let $A\in\mc{A}_1^*=\{A\in\sep:1\in A\}$, then $A$ is described uniquely by the gaps $a_1,a_2,\ldots,a_{r}\geq 1$ such that \[A=\{1,2+a_1,3+a_1+a_2,\ldots,r+a_1+\cdots+a_{r-1}\},\] and $\sum_{i=1}^{r}a_i=n-r$. So $|\A_1^*|$ is equal to the number of ways of choosing integers $b_1,b_2,\ldots,b_r \geq 0$ with $\sum_{i=1}^{r}b_i=n-2r$ which is simply
\[
\binom{n-2r+r-1}{r-1}=\binom{n-r-1}{r-1}.\]
\end{proof}
\begin{proof}[Proof of Lemma \ref{fequal}]
If $j\geq 2$ then $j\in f(A) \iff j+1 \in A$. Hence $f(A)=f(B)$ implies that
\[
A\cap\{3,4,\ldots,n\}=B\cap\{3,4,\ldots,n\}.\]
So if $A,B\in\sep$ and $f(A)=f(B)$ but $A\neq B$ then it must be the case that one of them contains $1$ and the other contains $2$.\end{proof}
\begin{proof}[Proof of Lemma \ref{fsep}]
Clearly $f(A)\in [n-1]^{(r)}$ (since nothing maps to $n$). So if $f(A)\not\in [n-1]_*^{(r)}$ then we must have either $1,2\in f(A)$ or $1,n-1 \in f(A)$. Hence either $1,3\in A$ or $2,n \in A$ and so $A\not\in \mc{B}\cup\mc{C}$. \end{proof}
\begin{proof}[Proof of Lemma \ref{sint}]
If $\mc{I}$ is an intersecting family and $A,B \in f(\mc{I})$ then there exist $C,D \in \mc{I}$ such that $A=f(C)$ and $B=f(D)$. Then $\emptyset\neq f(C\cap D) \subseteq f(C)\cap f(D)=A\cap B$. So $f(\mc{I})$ is intersecting.\end{proof}

\begin{proof}[Proof of Lemma \ref{slem}(a)]
Each set in $\mc{H}$ comes from applying $f$ to a set in $\sep$. The only way such a set can fail to belong to $[n-1]_*^{(r)}$ is if it contains $1$ but no set in $\mc{H}$ contains $1$. Hence $\mc{H}\subseteq [n-1]^{(r)}_*$.\end{proof}
\begin{proof}[Proof of Lemma \ref{slem}(b)]
The fact that the three families of sets $f(\mc{D})-\{1\}$, $f(\mc{E})-\{1\}$ and $\mc{F}-\{1\}$ are pairwise disjoint follows from considering how sets in these families meet the set $\{2,n-1\}$.

If $A\in f(\mc{D})-\{1\}$ then $2\in A$ and $n-1\not\in A$. If $A\in f(\mc{E})-\{1\}$ then $2\not\in A$ and $n-1\in A$. If $A\in \mc{F}-\{1\}$ then there exist $B\in\mc{B}$ and $C\in \mc{C}$ such that $f(B)=f(C)=A\cup\{1\}$. Then using Lemma \ref{fequal} we have $C=(B\b\{2\})\cup\{1\}$. So $n,3 \not\in B$, which implies that $2,n-1\not\in A$. Hence the three families are pairwise disjoint since for any $H\in \mc{H}$ we have
\[
H\cap\{2,n-1\}=\left\{\begin{array}{cc} \{2\},& H\in f(\mc{D})-\{1\}\\ \{n-1\},& H\in f(\mc{E})-\{1\} \\ \emptyset,& H\in \mc{F}-\{1\}.\end{array}\right.\]
\end{proof}
\begin{proof}[Proof of Lemma \ref{slem}(c)]
We wish to show that
\[
\mc{H}=(f(\mc{D})-\{1\})\cup (f(\mc{E})-\{1\})\cup(\mc{F}-\{1\}),\] is intersecting.

Let $A,B \in \mc{H}$, there are six cases to examine. We consider first the three cases when $A$ and $B$ belong to the same subfamily of $\mc{H}$. If $A,B\in f(\mc{D})-\{1\}$ then $2\in A\cap B$. Similarly if $A,B\in f(\mc{E})-\{1\}$ then $n-1\in A\cap B$. 

If $A,B \in \mc{F}-\{1\}$ then there are $B_1 \in \mc{B}$ and $C_1\in \mc{C}$ such that $A=f(C_1)\b\{1\}$ and $B=f(B_1)\b\{1\}$. Then $C_1,B_1\in \mc{A}\implies C_1\cap B_1\neq \emptyset$. Also $1\not\in B_1$ and $2\not\in C_1$ so there exists $j \in C_1\cap B_1$ with $j\geq 3$. Hence $f(j)\in (f(C_1)\b\{1\})\cap (f(B_1)\b\{1\})=A\cap B$. 

Now suppose $A\in f(\mc{D})-\{1\}$ and $B\in f(\mc{E})-\{1\}$. Then there exist $D\in\mc{D}$ and $E\in \mc{E}$ such that $A=f(D)\b\{1\}$ and $B=f(E)\b\{1\}$. Now $D,E \in \mc{A} \implies D\cap E \neq \emptyset$. Also  $2\not\in D$ and $1\not\in E$ so there exists $j\in D\cap E$ with $j \geq 3$. Hence $f(j)\in A\cap B$.

Next suppose $A\in f(\mc{D})-\{1\}$ and $B\in \mc{F}-\{1\}$. Then there exist $D\in \mc{D}$, $B_1\in\mc{B}$ such that $A=f(D)\b\{1\}$ and $B=f(B_1)\b\{1\}$. Then $B_1,D \in \mc{A}\implies B_1 \cap D\neq \emptyset$. Also $1\not\in B_1$ and $2\not\in D$ so there exists $j \in B_1\cap D$ with $j \geq 3$. Hence $f(j) \in A\cap B$.

Finally if $A\in f(\mc{E})-\{1\}$ and $B\in \mc{F}-\{1\}$, then there exist $E\in \mc{E}$ and $C_1\in \mc{C}$ such that $A=f(E)\b\{1\}$ and $B=f(C_1)\b\{1\}$. Since $C_1,E \in \mc{A}$ we have $C_1 \cap E\neq \emptyset$. Also $2\not\in C_1$ and $1\not\in E$ so there exists $j \in C_1\cap E$ with $j\geq 3$. Then $f(j) \in A\cap B$.

Hence $\mc{H}$ is an intersecting family. \end{proof}
\begin{proof}[Proof of Lemma \ref{slem}(d)]
We need to prove that $f(\mc{H})\subseteq [n-2]_*^{(r-1)}$.
We note first that by part (a) of this lemma we have $\mc{H}\subseteq [n-1]_*^{(r-1)}$. Let $H\in \mc{H}$ and consider $f(H)$. Clearly $f(H)\in [n-2]^{(r-1)}$. We simply need to check that we do not have $1,2 \in f(H)$ or $1,n-2 \in f(H)$. Now since $1\not\in H$ we have $1\in f(H)\iff 2\in H$. But we showed during the proof of part (b) of this lemma that $2\in H\implies H\in f(\mc{D})-\{1\}$. In which case there is $D\in \mc{D}$ such that $H=f(D)$. Now $D \in \mc{D} \implies 1,3 \in D \implies n,4\not\in D \implies n-1,3 \not\in f(D)=H \implies 2,n-2 \not\in f(H)$ as required. \end{proof}
\begin{proof}[Proof of Lemma \ref{uni1}]
Suppose that $f(\mc{H})\simeq \B_1^*=\{A\in [n-2]_*^{(r-1)}:1\in A\}$. We need to show that either $\A\simeq \A_1^*$ or there exists $j\in [n]$ such that $\A_{j,j+2}^*=\{A\in\sep:j,j+2\in A\}\subset \A$.

Since $f(\mc{H})\simeq \B_1^*$ there exists $i \in [n-2]$ such that $f(\mc{H})=\B_i^*=\{A\in [n-2]_*^{(r-1)}:i\in A\}$. If $i=1$ then every set in $\mc{H}$ contains $2$ and so $\mc{H}=f(\mc{D})-\{1\}$. Since $\mc{D}\subseteq \A_{1,3}^*$ and this last family has size $\binom{n-r-2}{r-2}=|f(\mc{H})|=|\mc{D}|$ then $\mc{D}=\A_{1,3}^*$. Similarly if $i=n-2$ then $\mc{E}=\A_{2,n}^*$. So if $i\in\{1,n-2\}$ then there is $j$ such that $\A_{j,j+2}^*\subset \A$. 

We will now show that if $i\in\{2,\ldots,n-3\}$ then $\A_{1,i+2}^*\cup\A_{2,i+2}^* \subseteq \A$. Suppose $A=\{1,a_2,\ldots,a_l,i+2,a_{l+2},\ldots,a_r\}\in \A_{1,i+2}^*$ but $A\not \in \A$. Then since $f(\mc{H})=\B_{i}^*$ it must be the case that $B=\{a_2-2,\ldots,a_l-2,i,a_{l+2}-2,\ldots,a_r-2\}\in f(\mc{H})$. Hence $C=(A\b\{1\})\cup\{2\}\in\A$. Now $1\in A$ implies that $n\not\in C$ and so $C\not\in\mc{D}\cup\mc{E}$. This implies that $B\in f(\mc{F})$ and so $A\in \A$, a contradiction. Hence $\A_{1,i+2}^*\subset \A$. Similarly $\A_{2,i+2}^*\subset \A$.

We can now show that every set in $\A$ must contain $i+2$ and so $\A\simeq \A_1^*$. Suppose not, then there exists a set $A=\{a_1,\ldots,a_r\}\in \A$ such that $i+2\not\in A$. But then choosing a single point from each gap of $A$ we may construct a set in $\A_{1,i+2}^*\cup\A_{2,i+2}^*\subseteq\A$, that is disjoint from $A$, contradicting the fact that $\A$ is intersecting. \end{proof}
\begin{proof}[Proof of Lemma \ref{uni2}]
Suppose equality holds in (\ref{bound}) and $\A_{1,3}^*\subseteq \A$. We show first that every set in $\A$ meets $\{1,3\}$. Suppose not, then there is a set $A=\{a_1,\ldots,a_r\}\in \A$ such that $A\cap\{1,3\}=\emptyset$. Then $a_3 \geq 6$ so $B=\{1,3,a_3-1,\ldots,a_r-1\}\in \A_{1,3}^*\subset \A$. But $A\cap B=\emptyset$ contradicting the fact that $\A$ is intersecting. Hence every set in $\A$ meets $\{1,3\}$. 

Let $\B_1=\{A\in\A:1\in A,3\not\in A\}$ and $\B_3=\{A\in\A:1\not\in A,3\in A\}$. For equality to hold in (\ref{bound}) we must have $|\B_1\cup\B_3|=\binom{n-r-2}{r-1}$. Let $\mc{C}_1=\{C\in\sep:1\in C,3\not\in C\}$ and $\mc{C}_3=\{C\in\sep:1\not\in C,3\in C\}$. For $A=\{1,a_2,\ldots,a_r\}\in \mc{C}_1$ define $g(A)=\{3,a_2+1,\ldots,a_r+1\}$. Clearly $A\cap g(A)=\emptyset$ and $g$ is a bijection from $\mc{C}_{1}$ to $\mc{C}_3$. Then, since $|\mc{C}_1|=|\mc{C}_3|=\binom{n-r-2}{r-1}$, it must be the case that for any $A\in \mc{C}_1$ exactly one of the sets $A$ and $g(A)$ belongs to $\A$. \end{proof}
\section{The Erd\H os-Ko-Rado Theorem for $k$-separated sets}
\begin{theorem}
\label{kCproof}
Let $k \geq 2$, $r \geq 1$, $n \geq (k+1)r$ and $\mc{A}\subseteq \sepk$ be intersecting then $|\A|\leq |\A_1^*|$ where $\mc{A}_1^*=\{A\in\sepk: 1 \in A \}$. Moreover $\A_1^*$ is the unique extremal family up to isomorphism.
\end{theorem}
\begin{proof}
We proceed by induction on $n$. First note that the result is trivial for $r\leq 2$ so we may suppose $r\geq 3$. Also we may suppose that $n \geq (k+1)r+1$ since if $n=(k+1)r$ then there are only $k+1$ sets in $\sepk$ and these are all pairwise disjoint so the result follows.

Let $\mc{A}\subseteq \sep$ be intersecting. Define the function $f:[n]\to [n-1]$ by 
\[
f(j)=\left\{\begin{array}{cc} 1,& j=1\\ j-1,& j\geq 2.\end{array}\right.
\]
Consider the following partition of $\A$
\[
\mc{A}=\mc{B}\cup\mc{C}\cup\bigcup_{i=0}^{k}\mc{D}_i,
\]
where
\begin{eqnarray*}
\mc{B}& =&\{A\in \A: 1\not\in A \tr{ and }f(A)\in [n-1]_{k}^{(r)}\},\\
\mc{C}& =&\{A\in \A: 1\in A \tr{ and }f(A)\in [n-1]_{k}^{(r)}\},\\
\mc{D}_0& =&\{A\in \A: 1,k+2\in A\},\\
\mc{D}_i&= &\{A\in \A:n+1-i,k+2-i\in A\},\ (1\leq i \leq k).
\end{eqnarray*}
Note that this is a partition of $\mc{A}$ since either $f(A)\in[n-1]_k^{(r)}$ and so $A\in \mc{B}\cup \mc{C}$ or $f(A)\not\in [n-1]_k^{(r)}$. If the latter holds then $1$ must lie in a gap of $A$ of size exactly $k$ or $1$ must be the left endpoint of such a gap and so there exists $0\leq i \leq k$ such that $A\in\mc{D}_i$.  

We will need to use the following results (again we defer their proofs until later).
\begin{lemma}
\label{kcval}
If $\mc{A}_1^*=\{A\in\sepk:1\in A\}$ then \[|\mc{A}_1^*|=\binom{n-kr-1}{r-1 }.\]
\end{lemma}
\begin{lemma}
\label{kfequal}
Suppose $A,B\in\sepk$ with $A\neq B$ and $f$ is as defined above. If $1\leq j\leq k$ and $f^j$ denotes $f$ iterated $j$ times then \[
f^j(A)=f^j(B) \implies A\Delta B=\{c,d\}\tr{ for }1\leq c < d \leq j+1.\]
\end{lemma}
Since no set in $\mc{B}$ contains $1$, Lemma \ref{kfequal} with $j=1$ implies that $|f(\mc{B})|=|\mc{B}|$. Similarly since no set in $\mc{C}$ contains $2$ we have $|f(\mc{C})|=|\mc{C}|$. Hence \[|\mc{B}|+|\mc{C}|=|f(\mc{B})|+|f(\mc{C})|=|f(\mc{B})\cup f(\mc{C})|+|f(\mc{B})\cap f(\mc{C})|.\]

Then, since $\mc{B}\cup\mc{C}$ is an intersecting family in $\sepk$, Lemma \ref{sint} implies that $f(\mc{B})\cup f(\mc{C})=f(\mc{B}\cup\mc{C})$ is an intersecting family in $[n-1]_k^{(r)}$. Hence by our inductive hypothesis for $n-1 \geq (k+1)r$ and Lemma \ref{kcval} we have 
\[
|f(\mc{B})\cup f(\mc{C})|\leq \binom{n-1-kr-1}{r-1}=\binom{n-kr-2}{r-1}.\]
Let\[
\mc{E}=f(\mc{B})\cap f(\mc{C}).\]
For any family of sets $\mc{G}$ recall that \[
\mc{G}-\{1\}=\{G\b\{1\}:G\in\mc{G}\}.\]
Define \[
\mc{F}=(f^{k-1}(\mc{E})-\{1\})\cup\bigcup_{i=0}^{k}(f^{k}(\mc{D}_i)-\{1\}),\]
where $f^j$ denotes the function $f$ iterated $j$ times.

We now require the following result.
\begin{lemma}
\label{kslem}
Let $\mc{D}_i,\mc{E},\mc{F}$ and $f$ be as defined above, then 
\begin{itemize}
\item[(a)]$\mc{F}\subseteq [n-k]^{(r-1)}_k$.
\item[(b)]$f^k(\mc{D}_0)-\{1\}$, $f^k(\mc{D}_1)-\{1\},\ldots,$ $f^k(\mc{D}_k)-\{1\}$ and $f^{k-1}(\mc{E})-\{1\}$ are pairwise disjoint families of sets.
\item[(c)] $\mc{F}$ is intersecting. 
\item[(d)] $f(\mc{F}) \subseteq [n-k-1]_k^{(r-1)}$. 
\end{itemize}
\end{lemma}
Since no set in $\mc{F}$ contains $1$, Lemma \ref{kfequal} with $j=1$ implies that $|\mc{F}|=|f(\mc{F})|$. Lemma \ref{kslem}(c) says that $\mc{F}$ is intersecting and so Lemma \ref{sint} implies that $f(\mc{F})$ is also intersecting. Then Lemma \ref{kslem}(d) tells us that $f(\mc{F}) \subseteq [n-k-1]_k^{(r-1)}$. So using Lemma \ref{kcval} and the inductive hypothesis for $n-k-1\geq (k+1)(r-1)$ we obtain\[
|\mc{F}|\leq \binom{n-k-1-k(r-1)-1}{(r-1)-1}=\binom{n-kr-2}{r-2}.\]
If $D\in\mc{D}_i$ then the definition of $\mc{D}_i$ implies that\[
 D\cap \{1,2,\ldots,k+1\}=\left\{\begin{array}{cc} \{1\}, & i=0\\\{k+2-i\},& i\geq 1.\end{array}\right.\]
Hence using Lemma \ref{kfequal} with $j=k$ we obtain $|\mc{D}_i|=|f^k(\mc{D}_i)|$ for $0 \leq i \leq k$.

If $E\in\mc{E}\subseteq f(\mc{C})\subseteq [n-1]_{k-1}^{(r)}$ then $E\cap \{1,2,\ldots,k\}=\{1\}$. So using Lemma \ref{kfequal} with $j=k-1$ we obtain $|\mc{E}|=|f^{k-1}(\mc{E})|$. Finally Lemma \ref{kslem}(b) tells us that
\begin{eqnarray*}
|\mc{F}|&=&|f^{k-1}(\mc{E})|+\sum_{i=0}^{k}|f^k(\mc{D}_i)|\\
 &=&|\mc{E}|+\sum_{i=0}^{k}|\mc{D}_i|.
\end{eqnarray*}
So we obtain
\begin{eqnarray}
\nonumber
|\mc{A}|&=&|\mc{B}|+|\mc{C}|+\sum_{i=0}^{k}|\mc{D}_i|\\
\nonumber
& = & |f(\mc{B})\cup f(\mc{C})|+|\mc{E}|+\sum_{i=0}^{k}|\mc{D}_i|\\
\nonumber
& = & |f(\mc{B})\cup f(\mc{C})|+|\mc{F}|\\
\label{kbound}
& \leq & \binom{n-kr-2}{r-1} + \binom{n-kr-2}{r-2}=\binom{n-kr-1}{r-1}=|\mc{A}_1^*|.
\end{eqnarray}
This completes the proof of the bound in Theorem \ref{kCproof}.

We now need to show that if equality holds in (\ref{kbound}) then $\A\simeq \A_1^*$. This follows in the same way as in the proof of Theorem \ref{Cproof} although the details are more involved. The only real difference is that for $k\geq 2$ the extremal family is always unique. The following two lemmas are obvious analogues of Lemmas \ref{uni1} and \ref{uni2}.
\begin{lemma}
\label{uni1k}
If equality holds in (\ref{kbound}) and the family $f(\mc{F})$ defined above is isomorphic to $\B_1^*=\{A\in [n-k-1]_k^{(r-1)}:1\in A\}$ then either $\A\simeq \A_1^*$ or there exists $j\in[n]$ such that $\A_{j,j+k+1}^*=\{A\in\sepk:j,j+k+1\in A \}\subset \A$.
\end{lemma}
\begin{lemma}
\label{uni2k}
If equality holds in (\ref{kbound}) and $A_{1,k+2}^*=\{A\in \sepk:1,k+2\in A\} \subset \A$ then every set in $\A$ meets $\{1,k+2\}$. Moreover if $A=\{1,a_2,\ldots,a_r\}\in\sepk$ satisfies $A\cap\{1,k+2\}=\{1\}$ then exactly one of the sets $A$ and $g(A)=\{k+2,a_2+k,\ldots,a_r+k\}$ belongs to $\A$. 
\end{lemma}
Now suppose equality holds in (\ref{kbound}). If $n=(k+1)r$ then clearly $\A\simeq \A_1^*$ so we may suppose that $n\geq (k+1)r+1$. We will prove that $\A\simeq \A_1^*$ by induction on $r$. Clearly the result holds for $r=1$ or $2$ so suppose $r\geq 3$. Since equality holds in (\ref{kbound}) we must have $|f(\mc{F})|=\binom{n-kr-2}{r-2}$ and so our inductive hypothesis for $r$ implies that $f(\mc{F})\simeq \B_1^*=\{A\in [n-k-1]_k^{(r-1)}:1\in A\}$. If $\A\not\simeq \A_1^*$ then Lemma \ref{uni1k} implies that there exists $j\in [n]$ such that $\A^*_{j,j+k+1}=\{A\in\sepk:j,j+k+1\in A\}\subset \A$. Without loss of generality we may suppose that $j=1$.

By the first part of Lemma \ref{uni2k} we know that every set in $\A$ meets $\{1,k+2\}$. Hence we may write $\mc{A}=\B_1\cup\B_{k+2}\cup\B_{1,k+2}$, where $\B_1=\{A\in\A:1\in A,k+2\not\in A\}$, $\B_{k+2}=\{A\in\A:1\not\in A,k+2\in A\}$ and $\B_{1,k+2}=\{A\in\A:1,k+2\in A\}$. For uniqueness to hold we need to show that one of the families $\B_1$ or $\B_{k+2}$ must be empty. We will prove that if both families are non-empty then $\A$ contains two disjoint sets, a contradiction.

Since we are supposing that equality holds in (\ref{kbound}), Lemma \ref{uni2k} tells us that if $A=\{1,a_2,\ldots,a_r\}\in\sepk$ and $\A\cap\{1,k+2\}=\{1\}$ then exactly one of the sets $A$ and $g(A)=\{k+2,a_2+k,\ldots,a_r+k\}$ belong to $\A$. We now wish to show that if $A\in \B_1$ and $B\in\sepk$ is obtained from $A$ by shifting one vertex of $A\b\{1\}$ one place clockwise then $B\in \B_1$. Suppose not, then by Lemma \ref{uni2k} we have $g(B)\in\B_{k+2}$. But $g(B)\cap A =\emptyset$, contradicting the fact that $\A$ is intersecting. Similarly if $A\in \B_{k+2}$ and $B\in\sepk$ is obtained from $A$ by shifting one vertex of $A\b \{k+2\}$ one place anti-clockwise then $B\in \B_{k+2}$.

Now suppose, for a contradiction, that both $\B_1$ and $\B_{k+2}$ are non-empty. Then, by shifting points one by one, we see that $A=\{k+2,2k+3,\ldots,(k+1)r+1\}$ and $B=\{1,n-(k+1)(r-1)+1,\ldots,n-2(k+1)+1,n-k\}$ both belong to $\A$. If $n\not\equiv 0\mod (k+1)$ then these are disjoint and we have a contradiction. So we may suppose that $n$ is a multiple of $(k+1)$ and hence $n\geq (k+1)(r+1)\geq (k+1)r+2$. Then either there is a set in $\B_1$ not containing $n-k$ and so $C=\{1,n-(r-1)(k+1),\ldots,n-2(k+1),n-(k+1)\}\in\A$ or $C\not\in\A$. If $C\in \A$ then we have a contradiction since $A\cap C=\emptyset$. Otherwise Lemma \ref{uni2k} implies that $g(C)=\{k+2,n-(r-2)(k+1)-1,n-(r-3)(k+1)-1,\ldots,n-1\}\in\A$. (Note that $C,g(C)\in\sepk$ follows from $n\geq (k+1)r+2$.) But since $k\geq 2$ we have $1\not\equiv -1 \mod (k+1)$ and so $B\cap g(C)=\emptyset$. This contradiction completes the proof of Theorem \ref{kCproof}.\end{proof}
We now prove the lemmas. The definitions of $\mc{A},\mc{B},\mc{C},\mc{D},\mc{E},\mc{F},\mc{H}$ and $f$ are as in the proof of Theorem \ref{kCproof}.

\begin{proof}[Proof of Lemma \ref{kcval}]
Let $A\in\mc{A}_1^*=\{A\in\sepk:1\in A\}$, then $A$ is described uniquely by the gaps $a_1,a_2,\ldots,a_{r}\geq k$ such that \[A=\{1,2+a_1,3+a_1+a_2,\ldots,r+a_1+\cdots+a_{r-1}\},\] and $\sum_{i=1}^{r}a_i=n-r$. So $|\A_1^*|$ is equal to the number of ways of choosing integers $b_1,b_2,\ldots,b_r \geq 0$ with $\sum_{i=1}^{r}b_i=n-(k+1)r$ which is simply
\[
\binom{n-(k+1)r +r -1}{r-1}=\binom{n-kr-1}{r-1}.\]
\end{proof}
\begin{proof}[Proof of Lemma \ref{kfequal}]
Let $A,B\in\sepk$ and suppose $f(A)=f(B)$ but $A\neq B$. If $a\geq 2$ then $a\in f^j(A)\iff a+j\in A$. Hence
\[
A\cap\{j+2,\ldots,n\}=B\cap\{j+2,\ldots,n\}.\]
So if $A\neq B$ then there exist $c,d\in\{1,\ldots,j+1\}$ such that $c\in A$ and $d\in B$ and $c\neq d$. But $j\leq k$ implies that $A\cap\{1,\ldots ,j+1\}=\{c\}$ and $B\cap\{1,\ldots,j+1\}=\{d\}$ or vice-versa and so $A\Delta B=\{c,d\}$ as required. \end{proof}
\begin{proof}[Proof of Lemma \ref{kslem}(a)] 
Let $G\in f^{k-1}(\mc{E})$. Then there exists $E\in \mc{E}$ such that $f^{k-1}(E)=G$. So there exist $B\in\mc{B}$ and $C\in\mc{C}$ such that $f(B)=f(C)=E$. Now Lemma \ref{kfequal} with $j=1$ implies that $B\Delta C=\{1,2\}$. Then $1\in C$ implies that $C\cap\{n-k+1,\ldots,n\}=\emptyset$. Hence $E\cap\{n-k,\ldots,n\}=\emptyset$. So $G\cap\{n-2k+1,\ldots, n\}=\emptyset$. Also $2\in B \implies E\cap \{2,\ldots,k+1\}=\emptyset$ and we have 
\begin{equation}
\label{e1}
G\cap(\{n-2k+1,\ldots,n-k\}\cup\{1,2\})=\{1\}.\end{equation}
 Hence $G\b\{1\}\in [n-k]_k^{(r-1)}$.

Now suppose $G\in f^{k}(\mc{D}_0)$. Then there exists $D\in\mc{D}_0$ such that $G=f^k(D)$ and $1,k+2\in D $. Hence $1 \in G$ and 
\begin{equation}
\label{e2}
G\cap(\{n-2k+1,\ldots,n-k\}\cup\{1,\ldots,k+2\})=\{1,2\}\end{equation}
So $G\b\{1\}\in[n-k]_k^{(r-1)}$.

Finally suppose that $1\leq i \leq k$ and $G\in f^{k}(\mc{D}_i)$. Then there is $D\in\mc{D}_i$ such that $G=f^k(D)$ and $n+1-i,k+2-i\in D $. Hence $1=f^k(k+2-i) \in G$ and 
\begin{equation}
\label{e3}
G\cap(\{n-(i+2k)+1,\ldots,n-k\}\cup\{1,\ldots,k+2-i\})=\{n+1-(i+k),1\}
\end{equation}
So again  $G\b\{1\}\in[n-k]_k^{(r-1)}$. Hence $\mc{F}\subseteq [n-k]_k^{(r-1)}$. \end{proof}
\begin{proof}[Proof of Lemma \ref{kslem}(b)]
The fact that these families of sets are pairwise disjoint follows from considering how sets from these families meet the set $\{n-2k+1,\ldots,n-k\}\cup\{2\}$.

First suppose $G\in f^{k-1}(\mc{E})-\{1\}$. Using (\ref{e1}) we have
\[G\cap(\{n-2k+1,\ldots,n-k\}\cup\{2\})=\emptyset.\]
Now suppose $G\in f^{k}(\mc{D}_0)-\{1\}$. Using (\ref{e2}) we have
\[G\cap(\{n-2k+1,\ldots,n-k\}\cup\{2\})=\{2\}.\]
Finally suppose $G\in f^{k}(\mc{D}_i)-\{1\}$ with $1\leq i \leq k$. Using (\ref{e3}) we have
\[G\cap(\{n-2k+1,\ldots,n-k\}\cup\{2\})=\{n+1-(i+k)\}.\]
Hence these families are pairwise disjoint.\end{proof}
\begin{proof}[Proof of Lemma \ref{kslem}(c)]
We wish to show that
\[
\mc{F}=(f^{k-1}(\mc{E})-\{1\})\cup\bigcup_{i=0}^{k}(f^{k}(\mc{D}_i)-\{1\}),\]
is intersecting. 

We will consider four cases. First suppose that $A,B\in f^k(\mc{D}_i)-\{1\}$. So there exist $D,E\in\mc{D}_i$ such that $f^k(D)\b\{1\}=A$ and $f^k(E)\b\{1\}=B$. If $i=0$ then $k+2\in D\cap E$ so $ 2 \in A\cap B$. Otherwise $1\leq i \leq k$ and $n+1-i\in D\cap E$ so $n+1-(i+k)\in A\cap B$.

Next suppose that $A,B\in f^{k-1}(\mc{E})-\{1\}$ then there exist $E,F\in \mc{E}$ such that $f^{k-1}(E)\b\{1\}=A$ and $f^{k-1}(F)\b\{1\}=B$. Now $E,F \in \mc{E}$ implies that there exist $B_1\in\mc{B}$ and $C_1 \in \mc{C}$ such that $f(B_1)=E$ and $f(C_1)=F$. Then $B_1,C_1\in \mc{A} \implies B_1\cap C_1\neq \emptyset$ and Lemma \ref{kfequal} implies that $2\in B_1,1\in C_1 $. So there exists $j\geq k+3$ such that $j\in B_1\cap C_1$. Hence $3\leq f^k(j)\in A\cap B$.

For the next case suppose $0 \leq i < j \leq k$, $A\in f^k(\mc{D}_i)-\{1\}$ and $B\in f^k(\mc{D}_j)-\{1\}$. Then there exist $C\in\mc{D}_i$ and $D\in\mc{D}_j$ such that $f^k(C)\b\{1\}=A$ and $f^k(D)\b\{1\}=B$. Now $j\geq 1 \implies D \cap \{1,\ldots,k+2\}=\{k+2-j\}$ and $k+2-j \not\in C$. But $C,D \in \mc{A} \implies C\cap D \neq \emptyset$. Hence there exists $l \geq k+3$ such that $l\in C\cap D$ and so $3\leq f^k(l)\in A\cap B$. 

For the last case suppose that $A\in f^{k-1}(\mc{E})-\{1\}$ and $B\in f^k(\mc{D}_i)-\{1\}$. Then there exist $D\in\mc{D}_i$ and $E\in \mc{E}$ such that $f^{k}(D)\b\{1\}=B$ and $f^{k-1}(E)\b\{1\}=A$. Furthermore there exist $B_1\in \mc{B}$ and $C_1\in \mc{C}$ such that $f(B_1)=f(C_1)=E$. Lemma \ref{kfequal} implies that $1\in C_1$ and $2\in B_1$. So $B_1\cap\{1,\ldots,k+2\}=\{2\}$ and $C_1\cap\{1,\ldots,k+1\}=\{1\}$. Also $D\in \mc{D}_i$ implies that \[
D\cap\{1,2,\ldots,k+2\}=\left\{\begin{array}{cc}\{1,k+2\},&i=0\\\{k+2-i\},& 1\leq i \leq k.\end{array}\right.\]
Now $B_1,C_1,D \in \mc{A}$ implies that $B_1\cap D$ and $C_1 \cap D$ are both nonempty. So if $i=0$ then there exists $j \geq k+3$ such that $j\in B_1\cap D$ and hence $3\leq f^k(j)\in A\cap B$. Otherwise $1\leq i \leq k$ and there exists $j \geq k+3$ such that $j\in C_1\cap D$ and hence $3\leq f^k(j)\in A\cap B$. 

Hence $\mc{F}$ is  intersecting. \end{proof}
\begin{proof}[Proof of Lemma \ref{kslem}(d)]
We need to prove that $f(\mc{F})\subseteq [n-k-1]_k^{(r-1)}$.
We note first that by part (a) of this lemma we have $\mc{F}\subseteq [n-k]_k^{(r-1)}$. Let $F\in \mc{F}$ and consider $f(F)$. Clearly $f(F)\in [n-k-1]^{(r-1)}$. We simply need to check that $f(F)$  does not contain a gap of size exactly $k-1$ (as a subset of $[n-k-1]$). This will follow if we show that $F$ does not contain a gap of size exactly $k$ (as a subset of $[n-k]$) around $1$. More precisely we need to check that $F$ does not contain any of the following pairs of points:\[
(n-2k+1,2),\ldots,(n-k,k+1),(1,k+2).\]
Since $1\not\in F$ we know that the last pair in this list cannot belong to $F$. Also using (\ref{e1}) and (\ref{e2}) we know that if $F\in f^{k-1}(\mc{E})-\{1\}$ or $F\in f^{k}(\mc{D}_0)-\{1\}$ then 
\[
F\cap\{n-2k+1,\ldots,n-k\}=\emptyset.\]
Finally using (\ref{e3}) we know that if $F\in f^{k}(\mc{D}_i)-\{1\}$, with $1\leq i \leq k$, then
\[
F\cap(\{n-2k+1,\ldots,n-k\}\cup\{1,\ldots,k+2-i\})= \{n+1-(i+k)\}.\] 
Hence if $n-2k+j \in F\in f^k(\mc{D}_i)-\{1\}$, with $1\leq j \leq k$, then $j=k+1-i$ and so $j+1=k+2-i \not\in F$. Hence $F$ cannot contain any of the pairs of points in the above list (since they are all of the form $(n-2k+j,j+1)$). So $f(\mc{F})\subseteq [n-k-1]_k^{(r-1)}$ as required.\end{proof}
\begin{proof}[Proof of Lemma \ref{uni1k}]
Suppose that $f(\mc{F})\simeq \B_1^*=\{A\in [n-k-1]_k^{(r-1)}:1\in A\}$. We need to show that either $\A\simeq \A_1^*$ or there exists $j\in [n]$ such that $\A_{j,j+k+1}^*=\{A\in\sepk:j,j+k+1\in A\}\subset \A$.

Since $f(\mc{F})\simeq \B_1^*$ there exists $i \in [n-k-1]$ such that $f(\mc{F})=\B_i^*=\{A\in [n-k-1]_k^{(r-1)}:i\in A\}$. If $i=1$ then every set in $\mc{F}$ contains $2$ and so $\mc{F}=f^k(\mc{D}_0)-\{1\}$. Since $\mc{D}_0\subseteq \A_{1,k+2}^*=\{A\in\sepk:1,k+2\in A\}$ and this last family has size $\binom{n-kr-2}{r-2}=|f(\mc{F})|=|\mc{D}_0|$ then $\mc{D}_0=\A_{1,k+2}^*$. Similarly if $i\in\{n-2k,\ldots,n-k-1\}$ then $\A_{j,j+k+1}^*\subset \A$, where $j=i+k+1$. 

So suppose now that $i\in\{2,\ldots,n-2k-1\}$. We will first show that $\A_{1,i+k+1}^*\cup\A_{2,i+k+1}^* \subseteq \A$. Suppose $A=\{1,a_2,\ldots,a_r\}\in \A_{1,i+k+1}^*$ but $A\not \in \A$. Then since $f(\mc{F})=\B_{i}^*$ it must be the case that $B=\{a_2-(k+1),\ldots,a_r-(k+1)\}\in f(\mc{F})$. Hence $C=(A\b\{1\})\cup\{j\}\in\A$ for some $j\in\{2,\ldots,k+1\}$ and then $1\in A$ implies that $C\cap\{n-k+1,\ldots,n\}=\emptyset$. Hence $C\not\in\bigcup\limits_{i=0}^{k}\mc{D}_i$. This implies that $B\in f(f^{k-1}(\mc{E})-\{1\})$ and so $A\in \A$, a contradiction. Hence $\A_{1,i+k+1}^*\subset \A$. Similarly $\A_{2,i+k+1}^*\subset \A$.

Define $\mc{C}^*_j=\{A\in \sepk:i+k+1,j+1,n-k+j\in A\}$. So $\mc{D}_{j}\subseteq \mc{C}_{k+1-j}$. We claim that if $1\leq j\leq k+1$ then $\mc{C}_j^*\subset \A$. Fix $j$ and let $C\in\mc{C}^*_j$. Consider $B=f(f^k(C)\b\{1\})\in \B^*_i$. Then $n-2k+j\in B$ implies that $C\in\mc{D}_{j-(k+1)}\subseteq \A$.

To summarise, we now know that
 \begin{equation}
\label{ina}
\A_{1,i+k+1}^*\cup\A_{2,i+k+1}^*\cup\bigcup\limits_{j=1}^{k+1}\mc{C}_{j}^* \subseteq \A.
\end{equation}
We can now show that every set in $\A$ must contain $i+k+1$ and so $\A\simeq \A_1^*$. Suppose not, then there exists a set $A=\{a_1,\ldots,a_r\}\in \A$ such that $i+k+1\not\in A$. We will construct a set $B\in\A$ such that $A\cap B=\emptyset$. 

Suppose $a_j<i+k+1<a_{j+1}$. We consider first the case that $a_j=a_r$ and $a_{j+1}=a_1$. If $1\leq i+k+1 <a_1$ then $i+k+1\geq k+3$. So let $B=\{1,i+k+1,a_2-1,\ldots,a_{r-1}-1\}\in \A^*_{1,i+k+1}\subset \A$. Then $A\cap B = \emptyset$. Otherwise we have $i+k+1 \leq n$ and hence $i+k+1 \leq n-k$. So let $B=\{2,a_2+1,\ldots,a_{r-1}+1,i+k+1\}$ or $B=\{1,a_2+1,\ldots,a_{r-1}+1,i+k+1\}$ depending on whether or not $a_1=1$. In either case $B\in \A^*_{1,i+k+1}\cup\A^*_{2,i+k+1}\subset \A$ and $A\cap B = \emptyset$.

We turn to the general case, i.e. $a_j<i+k+1<a_{j+1}$ for $1\leq j \leq r-1$. Let the gap between $a_r$ and $a_1$ be of size $\alpha\geq k$, and suppose $a_j+\beta=a_{j+1}-\gamma=i+k+1$. We must consider two sub-cases. First suppose that $\max\{\beta,\gamma\}\leq k$. Let $C=\{a_1+\beta,\ldots a_{j-1}+\beta,i+k+1,a_{j+2}-\gamma,\ldots,a_r-\gamma\}$. So $i+k+1\in C$ and $A \cap C = \emptyset$. The gap between the first and last elements of $C$ is $\alpha+\beta+\gamma \geq 2k+1$. If $a_r-\gamma+k+1 \geq 1$ let $B=C\cup\{a_r-\gamma+k+1\}$. Otherwise if $a_1+\beta-(k+1)\leq 2$ let $B=C\cup\{a_1+\beta-(k+1)\}$. If neither holds then let $B=C\cup\{1\}$. In each case $A\cap B = \emptyset$ and (\ref{ina}) implies that $B \in \A$ as required.

For the second sub-case we suppose $\max\{\beta,\gamma\}\geq k+1$. Without loss of generality assume $\beta \geq k+1$. Let $\delta=\max\{k+1-\gamma,1\}$. We now define $C=\{a_1-1,\ldots,a_j-1,i+k+1,a_{j+1}+\delta,\ldots,a_{r-1}+\delta\}$. Now either $a_1=1$ in which case let $B=(C\b\{a_2-1\})\cup\{k+1\}$ or $a_1\geq 2$. If $a_1=2$ or $a_1=3$ then let $B=C\in\A_{1,i+k+1}^*\cup\A_{2,i+k+1}^* \subseteq \A$ and $A\cap B=\emptyset$. So suppose $a_1\geq 4$. In this case we construct $B$ from $C$ by shifting the last element of $C$ clockwise and/or the first element of $C$ anticlockwise until we obtain a set in $\cup_{j=1}^{k+1}\mc{C}^*_j$ that is disjoint from $A$. Again (\ref{ina}) implies that $B\in\A$ as required. 

Hence every set in $\A$ contains $i+k+1$ and so $\A\simeq \A_1^*$. \end{proof}
\begin{proof}[Proof of Lemma \ref{uni2k}]
Suppose equality holds in (\ref{kbound}) and $\A_{1,k+2}^*\subseteq \A$. We show first that every set in $\A$ meets $\{1,k+2\}$. Suppose not, then there is a set $A=\{a_1,\ldots,a_r\}\in \A$ such that $A\cap\{1,k+2\}=\emptyset$. Then $a_3 \geq 2k+4$ so $B=\{1,k+2,a_3-1,\ldots,a_r-1\}\in \A_{1,k+2}^*\subset \A$. But $A\cap B=\emptyset$ contradicting the fact that $\A$ is intersecting. Hence every set in $\A$ meets $\{1,k+2\}$. 

Let $\B_1=\{A\in\A:1\in A,k+2\not\in A\}$ and $\B_{k+2}=\{A\in\A:1\not\in A,k+2\in A\}$. Then for equality to hold in (\ref{kbound}) we must have $|\B_1\cup\B_{k+2}|=\binom{n-kr-2}{r-1}$. Let $\mc{C}_1=\{C\in\sepk:1\in C,k+2\not\in C\}$ and $\mc{C}_{k+2}=\{C\in\sepk:1\not\in C,k+2\in C\}$. Then for $A=\{1,a_2,\ldots,a_r\}\in \mc{C}_1$ define $g(A)=\{k+2,a_2+k,\ldots,a_r+k\}$. Clearly $A\cap g(A)=\emptyset$ and $g$ is a bijection from $\mc{C}_{1}$ to $\mc{C}_{k+2}$. Then since $|\mc{C}_1|=|\mc{C}_{k+2}|=\binom{n-kr-2}{r-1}$ it must be the case that for any $A\in \mc{C}_1$ exactly one of the sets $A$ and $g(A)$ belongs to $\A$. \end{proof}
\section{The Erd\H os-Ko-Rado theorem for weighted $k$-separated sets}
We conclude this paper with a short result for suitably weighted $k$-separated sets. Unlike our other results this follows simply from the original Erd\H os-Ko-Rado theorem (Theorem \ref{EKRthm}).

For $A=\{a_{1}, \dots , a_{r}\} \in [n]^{(r)}_{k}$ we define the \emph{weight} of $A$ to be
\[
w(A)=\prod_{i=1}^{r}\binom{a_{i+1}-a_{i}-1}{k},
\]
where $a_{r+1}=a_{1}+n$.
So the weight of a $k$-separated set $A\in \sepk$ is simply the number of different ways $A$ may be extended to form a set $B \in [n]^{((k+1)r)}$ by inserting exactly $k$ new elements into each gap in $A$. We then define the \emph{weight} of a family of sets $\mathcal{A} \subseteq [n]^{(r)}_{k}$ to be
\[
w(\mathcal{A})=\sum_{A \in \mathcal{A}}w(A).
\]
The following result says that an analogue of the Erd\H os-Ko-Rado theorem holds for weighted $k$-separated sets.
\begin{theorem}
\label{Wthm}
Let $n \geq 2(k+1)r$. If $\A \subseteq \sepk$ is intersecting then $w(\A)\leq w(\A^{*}_{1})$, where $\A_{1}^{*}=\{A\in \sepk : 1\in A\}$.
\end{theorem}
\begin{proof}
Consider the bipartite graph $G=(V\cup W,E)$ with vertex classes $V=\sepk$ and $W=[n]^{((k+1)r)}$. We define $E$ as follows: let $A \in V$ be adjacent to $B \in W$ if we can construct $B$ from $A$ by inserting exactly $k$ elements into each gap in $A$. 

Let $\mathcal{A}\subseteq \sepk$ be intersecting then 
\[
\Gamma(\mathcal{A})=\{B\in [n]^{((k+1)r)}:\exists A\in\A\tr{ such that }(A,B)\in E\}
\]
is also intersecting. Then, since $\Gamma(\mc{A})$ is an intersecting family of $(k+1)r$-sets from $[n]^{((k+1)r)}$ and $n \geq 2(k+1)r$, Theorem \ref{EKRthm} implies that  
\[
|\Gamma(\mathcal{A})| \leq \binom{n-1}{(k+1)r-1}.
\]
For distinct $A_{1},A_{2} \in \mathcal{A}$ we have $\Gamma(A_{1}) \cap \Gamma(A_{2}) = \emptyset$. To see this, suppose we had $B \in \Gamma(A_{1}) \cap \Gamma(A_{2})$ with $B=\{b_{1},\dots ,b_{(k+1)r} \}$. Without loss of generality we may suppose that $A_{1}=\{b_{1},b_{k+2}, \dots ,b_{(k+1)r-k}\}$ and $A_{2}=\{b_{i},b_{k+i+1}, \dots,$ $b_{(k+1)(r-1)+i}\}$, for some $2 \leq i \leq k+1$. Hence $A_{1} \cap A_{2} = \emptyset$. This contradicts the fact that $\mathcal{A}$ is intersecting. 

Then, since $|\Gamma(A)|=w(A)$, we have
\[
w(\mathcal{A})=\sum_{A \in \mathcal{A}}w(A)=\sum_{A \in \mathcal{A}} |\Gamma(A)|=|\Gamma(\mathcal{A})| \leq \binom{n-1}{(k+1)r-1}=w(\mathcal{A}_{1}^{\ast}).
\]
\end{proof}

\end{document}